\documentclass{ws-sd}
\usepackage{amssymb}
\usepackage{color, amsmath,amssymb, amsfonts, amstext, latexsym}

\usepackage{amssymb, epsfig, amssymb, latexsym}
\usepackage{amsmath}
\usepackage{graphicx}
\usepackage{longtable}

\newcommand{\Om}{\Omega}
\renewcommand{\phi}{\varphi}

\newcommand{\R}{{\mathbb R}}

   \newcommand{\PX}{{\Bbb{P}}}

\newcommand{\cF}{{\cal F}}

\newcommand{\hse}{\ensuremath{l^s(\xi,\Phi^0,\omega)}}
\newcommand{\hsk}[1]{\ensuremath{l^{(#1)}(\xi,\Phi^0,\omega)}}

\newcommand{\sz}{\ensuremath{ {\int_s^t z(\theta_r (\omega))\,dr}}}
\newcommand{\sZ}{\ensuremath{ {\int_s^t Z(\theta_r (\omega))\,dr}}}
\newcommand{\zz}[1]{\ensuremath{{z(\theta_{#1} (\omega))}}}
\newcommand{\ZZ}[1]{\ensuremath{{Z(\theta_{#1} (\omega))}}}
\newcommand{\fu}[1]{\ensuremath{{F_\psi^{\psi^{(d)} (#1)}}}}
\newcommand{\fus}[2]{\ensuremath{{\int^0_{#2} F_\psi^{\psi^{(d)} (#1)}\,d{#1}}}}
\newcommand{\fuss}[2]{\ensuremath{{\int^{#2}_0 F_\psi^{\psi^{(d)} (#1)}\,d{#1}}}}

\newcommand{\rb}{\right)}
\newcommand{\lb}{\left(}
\newcommand{\rB}{\right]}
\newcommand{\lB}{\left[}

\newcommand{\BE}{\begin{equation}}
\newcommand{\EE}{\end{equation}}
\newcommand{\BEN}{\begin{equation*}}
\newcommand{\EEN}{\end{equation*}}
\newcommand{\BAL}{\begin{align}}
\newcommand{\EAL}{\end{align}}
\newcommand{\BAN}{\begin{align*}}
\newcommand{\EAN}{\end{align*}}

\begin{document}
\markboth{X. Sun, X. Kan and J. Duan}{Approximation of invariant foliations}

\title{Approximation of invariant foliations\\ for stochastic   dynamical systems
\footnote{Part of this work was done while J. Duan was participating the Stochastic Partial Differential Equations programme at the Isaac Newton Institute for Mathematical Sciences, Cambridge, UK.  This work was partly supported by   NSF of China grants
  10971225  and  11028102, the NSF Grants  1025422 and
0731201,    the Cheung Kong Scholars Program, and an open research grant from the State Key Laboratory for Nonlinear Mechanics at the Chinese Academy of Sciences.
  }
}

\author{Xu Sun}

\address{School of Mathematics and Statistics\\Huazhong University of Science and Technology \\
  Wuhan 430074, China\\
xsun15@gmail.com}

\author{Xingye Kan}

\address{Department of Applied Mathematics\\ Illinois Institute of Technology \\
  Chicago, IL 60616, USA\\
xkan@iit.edu}

\author{Jinqiao Duan}

\address{Department of Applied Mathematics\\ Illinois Institute of Technology \\
  Chicago, IL 60616, USA\\
duan@iit.edu}

\maketitle

\begin{abstract}
Invariant foliations are geometric structures for describing and understanding  the qualitative behaviors of nonlinear dynamical systems. For stochastic dynamical systems, however,
these geometric structures themselves  are complicated  random sets. Thus it is desirable to have some techniques to approximate random invariant foliations.   In this paper, invariant foliations are approximated for   dynamical systems with small noisy perturbations, via asymptotic analysis. Namely,
   random invariant foliations are represented as a perturbation of the deterministic invariant foliations, with deviation errors estimated.

\end{abstract}


 \bigskip

%\emph{Dedicated to Professor Peter Imkeller on the occasion  of his 60th birthday }

\section {Introduction and motivation}
\quad Invariant foliations, as well as   invariant manifolds, provide geometric structures for understanding the qualitative behaviors of nonlinear
dynamical systems, and they have been extensively studied for deterministic systems \cite{Hirsch,GH,Bates}.

Invariant manifolds or  foliations for finite dimensional stochastic systems or stochastic differential equations (SDEs) were studied in
\cite{Wanner,Arnold,Mohammed,DuDuan}. Recently, the existence of  invariant manifolds  and invariant foliations for stochastic partial differential
equations was investigated in \cite{DuanLuSchm,DuanLuSchm2,Mohammed2,CDLS} and \cite{luschmfuss2008}, respectively. In \cite{XuSun}, we estimated the
impact of small noise on invariant manifolds for nonlinear systems. Note that random center-like invariant manifolds were approximated for some
stochastic differential equations (SPDEs) by Wang and Duan \cite{WangDuan}, and  Blomker and Wang \cite{BW}. In this paper, we consider a procedure
to approximate   invariant foliations for nonlinear systems perturbed by small noise. We   compare invariant foliations for the original
deterministic dynamical systems and for the randomly perturbed systems.

We consider the following nonlinear stochastic evolutionary equation with a multiplicative noise, in a separable Hilbert space $H$ with a scalar
product $<\cdot, \cdot>$ and the induced norm $\| \cdot \|= \sqrt{<\cdot, \cdot>}$:
\begin{align}\label{spde}
\frac{dU}{dt}=AU+F(U)+\epsilon \; U \circ  \dot{W},
\end{align}
where $A$ is a linear (bounded or unbounded) operator, ``$\circ$" is in the sense of Stratonovich stochastic calculus, $W=W(t, \omega)$ is a scalar Brownian motion defined
on a probability space $(\Om, \cF, \PX)$,  and $\epsilon$ is a positive parameter representing the intensity of the noise.
This covers some SDEs and SPDEs.  Note that the Ito's form
of (\ref{spde}) is
\begin{align*}
\frac{dU}{dt}=AU+\epsilon \frac{U}{2} +F(U)+\epsilon \; U \dot{W}.
\end{align*}
 The nonlinearity $F(U)$
satisfies $F(0)=0$ and $DF(0)=0$ and is Lipschitz continuous on $H$
\[
\|F(U_1)-F(U_2)\| \le L_F \|U_1-U_2\|,
\]
where $L_F$ is the positive Lipschitz constant and $\|\cdot\|$ the norm in the Hilbert space $H$. When the  nonlinearity $F(U)$ is locally Lipschitz continuous, the
approximation result in this paper can be applied to the modified stochastic equation where the nonlinearity is appropriately cut-off and thus obtain
approximation information for the local random invariant foliations.  The state space $H$ is the Euclidean space $\R^n$ when the above equation is a
SDE or a function space if the above equation is a SPDE. When $\epsilon=0$, Eq.~(\ref{spde}) reduces to a deterministic evolutionary equation:
\begin{align}\label{pde}
\frac{dU}{dt}=AU+F(U).
\end{align}
We compare the invariant foliations for the original deterministic system \eqref{pde} and for the randomly perturbed system \eqref{spde}, and
quantify their difference when the noise intensity $\epsilon$ is small.

This paper is organized as follows. In section 2, we review some basic concepts of random dynamical systems, and recall the existence result for
random invariant foliations. The main result on asymptotic analysis for random invariant foliations   is described in section 3, and two illustrative
examples are presented in section 4.

\section{Random invariant foliation}

 Following \cite{DuanLuSchm2}, we assume throughout the paper that the linear operator $A: D(A)\rightarrow H$
generates a strongly continuous semigroup $e^{At}$ on $H$, which satisfies the pseudo exponential dichotomy with exponents $\alpha >0>\beta$ and a
bound $K>0$, i.e., there exists a continuous projection $P^u$ on $H$ such that
\begin{enumerate}
\item[(i)] $P^ue^{At}=e^{At}P^u$;\\
\item[(ii)] The restriction $e^{At}|_{R(P^u)}$,  $t\ge 0$, is an isomorphism of the range $R(P^u)$ of $P^u$ onto itself, and we define $e^{At}$ for $t<0$ as
the inverse map;
\item[(iii)] The following estimates hold
\begin{align}\label{cdt}
\mid e^{At} P^u x\mid &\le Ke^{\alpha t} \mid x\mid, &&t\le 0,\\
\mid e^{At} P^s x\mid &\le Ke^{\beta t} \mid x\mid,  &&t\ge 0.\notag
\end{align}
where $P^s=I-P^u$. Denote $H^s=P^s H$, $H^u=P^u H$ and hence $H=H^s \bigoplus H^u$.
\end{enumerate}

\subsection{Random dynamical systems}
A measurable random dynamical system on Hilbert space $(H,\mathcal{B})$ over a driving system  $ (\theta(t))_{t\in T}$ with time $T$ is a mapping
\begin{align}
\nonumber \phi: T\times \Omega \times H\rightarrow H, (t,\omega,x)\rightarrow \phi (t,\omega,x),
\end{align}
with the following properties \cite{Arnold}:
\begin{enumerate}
\item[(i)] Measurability: $\phi$ is $\mathcal {B}(T)\otimes \mathcal{F}\otimes\mathcal{B}$--measurable.
\item[(ii)] The mappings $\phi(t,\omega)=\phi(t,\omega,\cdot): H\to H$ form a cocycle over $\theta(\cdot)$, i.e. they satisfy $\phi(0,\omega)=id_X$ for all
$\omega\in\Omega$ and $\varphi(t+s,\omega)=\varphi(t,\theta(s)\omega) \; \varphi(s,\omega)$ for all $s, t \in T$ and $\omega\in\Omega$.
\end{enumerate}
For SDEs and SPDEs \cite{Arnold}, we   identify $\omega(t)=W(t,\omega)$, and define the driving system   $\theta(t)$ is the  Wiener shift, i.e.,
$\theta_t\omega(\tau)=\omega(\tau+t)-\omega(t)$.

To facilitate   random dynamical systems study of  (\ref{spde}), we convert  it into an evolutionary equation with random coefficients, called
a random evolutionary equation. To this end,
 we introduce $z(w)$ as the stationary solution of the following Langevin   equation
\begin{align}
dz+zdt=\epsilon dW.   \nonumber
\end{align}
Then $z(w)=\epsilon Z(\omega)$, where   $Z(\omega)$ is the stationary solution of $dZ(t)+Z(t)dt=dW(t)$ and it can be expressed as
\[
 Z(\omega) = \int^0_{-\infty} e^\tau  \,d W(\tau).
\]
Moreover,
\[
 Z(\theta_t \omega ) =e^{-t}Z(\omega) +e^{-t} \int_0^t e^\tau   \,d W(\tau).
\]
 Define a transform
\begin{align}
\bar x: = T(\omega,\bar X)=\bar Xe^{-z(\omega)}\nonumber
\end{align}
with its inverse transform
\begin{align}\label{trnsfm}
\bar X: = T^{-1}(\omega,\bar x)=\bar x e^{z(\omega)}.
\end{align}
Denote $U(t,\omega,X)$ as the solution of (\ref{spde}) with initial value $X$. Introducing
\[
u=T(\theta_t \omega, U(t,\omega,X)=e^{-z(\theta_t \omega)}U(t,\omega,X),
\]
then the new system state $u$ satisfies the following random   evolutionary equation \cite{DuanLuSchm2}
\begin{align}\label{rde}
\frac{du}{dt}=Au+z(\theta_t \omega)u+G(\theta_t\omega,u),\;\;\; u(0)=x\in H,
\end{align}
where
\[
x=T(\omega,X)=e^{-z(\omega)}X
\]
and
\begin{align}\label{E.dfnG}
G(\omega,u):=e^{-z(\omega)}F(e^{ z(\omega)}u).
\end{align}
We often denote the  solution of \eqref{rde} to be $u= \phi(t, \omega, x)$.
The solution mapping of (\ref{rde}), i.e. $(t,\omega,x)\rightarrow  \phi(t,\omega,x)$, generates a random dynamical system. Thus (see
\cite{DuanLuSchm2})
\[
(t,\omega,X) \to T^{-1} (\theta_t \omega,  \phi (t,\omega,T(\omega,X)):= U(t,\omega,X)
\]
is also a random dynamical system. In fact, The relationship between solutions of (\ref{spde}) and (\ref{rde}) is described by
\begin{align*}
 U(t,\omega,X)&=T^{-1} (\theta_t \omega, \phi (t,\omega,T(\omega,X))\\
 \phi(t,\omega,x)&=T      (\theta_t \omega, U(t,\omega,T^{-1}(\omega,x)).
 \end{align*}

\subsection{Random invariant foliation}

The concept of invariant foliation is about quantifying certain sets (called leaves or fibers) in state space $H$, starting from all points in such a leaf the dynamical orbits have similar asymptotic behaviors. These leaves are thus building blocks for understanding dynamics.

Let us consider a leaf for random invariant foliation for the above random dynamical system $\phi (t, \omega, x)$.
\quad A   leaf passing through a point $\Phi^0$ in the state space $H$,   denoted as $W(\Phi^0, \omega)$, is a random set and is invariant  in the following special sense
  \cite{Arnold,luschmfuss2008}
\begin{align}
\phi(t,\omega,W(\Phi^0, \omega))\subset W(\phi(t, \omega, \Phi^0), \theta_t \omega)\;\;for\;\;t\ge 0.  \nonumber
\end{align}
If we can represent $Q$ as a graph of a $C^k$ (or Lipschitz) mapping, then $Q(\Phi^0, \omega)$ is called a $C^k$ (or Lipschitz) leaf for the random invariant foliation.
The existence of random invariant foliation for \eqref{rde} is shown in  \cite{luschmfuss2008}. To facilitate our asymptotic analysis in the next section, we recall as follows.
We only consider stable leaves, still denoted as $W(\Phi^0, \omega)$. Unstable leaves may be considered similarly.

Define
\begin{align}\label{df_psi}
\psi(t)=\tilde \Phi(t)-\Phi(t),
\end{align}
where $\tilde \Phi(t)=\phi(t, \omega, \tilde \Phi^0)$ and $\Phi(t)=\phi(t, \omega,  \Phi^0)$ are two solutions of \eqref{rde} starting at two initial states
$\tilde \Phi^0$ and $ \Phi^0$, respectively. Also introduce   the following Banach space, for each $\eta$,
$\beta<\eta<\alpha$,
\begin{align}
\hat C^{+}_{\eta}= \{\phi :  [0,\infty )\rightarrow H\mid \phi \;\;is\;\;continuous\;\;and\;\;\underset{t\in  [0,\infty )}{sup} \;e^{-\eta t
-\int_0^t z(\theta_{\tau}\omega)d\tau}\| \phi (t)\| <\infty\}\nonumber
\end{align}
with the norm
\begin{align*}
\| \phi\|_{\hat C^{+}_{\eta}}=\underset{t\in  [0,\infty )}{sup} \;e^{-\eta t -\int_0^t z(\theta_{\tau}\omega)d\tau}\| \phi (t)\|.
\end{align*}
It is shown in (\cite{luschmfuss2008}) that $\tilde \Phi^0 \in W(\Phi^0, \omega)$ if and only if there exists a function $\psi(\cdot) \in \hat
C^{+}_{\eta}$ with $\psi(0)=\tilde \Phi^0 - \Phi^0$ and
\begin{align}\label{exp_u1tmp}
\psi(t)&=&e^{At+\int_0^t z(\theta_s\omega)ds}\xi+\int_0^t e^{A(t-s)+\int_s^t z(\theta_r\omega)dr}P^s \Delta G(\theta_s\omega, \psi(t), \Phi(t))ds\notag\\
&&+\int_\infty^t e^{A(t-s)+\int_s^t z(\theta_r\omega)dr}P^u  \Delta G(\theta_s\omega,\psi(t), \Phi(t))ds
\end{align}
where $\xi=P^s(\tilde \Phi^0-\Phi^0)$, and $\Delta G(\omega, \psi, \Phi)=G(\omega, \psi+\Phi)- G(\omega, \Phi)$. Under the gap condition
\begin{align*}
K=CL\lb \frac{1}{\alpha - \eta}+\frac{1}{\eta-\beta}\rb <1,
\end{align*}
  there exists an invariant foliation for (\ref{rde}) whose stable leaf is given by
\begin{align*}
W (\Phi^0, \omega)=\left\{ \xi + l^s (\xi, \Phi^0, \omega) \mid \xi\in H^s \right\},
\end{align*}
where $\Phi^0 \in H$, $(\xi,\Phi^0,\omega)\to l^s (\xi, \Phi^0, \omega)$ is measurable and Lipschitz continuous in $\xi$ and
\begin{align*}
\phi(t,\omega,W(\Phi^0,\omega)) \subset W (\Phi(t,\phi^0,\omega),\theta_t \omega).
\end{align*}
Moreover
\begin{align}\label{star2}
l^s (\xi, \Phi^0, \omega)=P^u \Phi^0 + P^u \psi(0;\xi -P^s \Phi^0, \Phi^0, \omega), \quad \xi \in H^s,
\end{align}
where
\begin{align*}
P^u \psi(t,\xi,\Phi^0,\omega)& = \int_\infty^t e^{A(t-s) + \int_s^t z(\theta_r\omega)\,dr} P^u (G(\theta_s \omega, \psi(s,\xi,\Phi^0, \omega) +
\Phi(s)- G(\theta_s \omega, \Phi(s))\,ds.
\end{align*}

%Given a leaf $W(\omega)$ for the invariant foliation for \eqref{rde}, the corresponding leaf of the invariant foliation for the original equation (\ref{spde}) can be expressed as %\cite{DuanLuSchm2}
%\begin{align}\label{E.trans}
%\tilde{W}(\omega) =T^{-1}(\omega,W(\omega)),
%\end{align}
%with $T^{-1}$ being defined   in (\ref{trnsfm}).
It also follows from (\ref{exp_u1tmp}) that $\psi(t)$, as defined in \eqref{df_psi}, satisfies the following equation
\begin{align}\label{star_1}
\frac{d\psi}{dt}=Au+z(\theta_t \omega)\psi + \Delta G(\theta_t\omega,\psi(t), \Phi(t)).
\end{align}

\section{Asymptotic analysis for random invariant foliation}

In this section, we   propose an approach to approximate the random invariant foliation by asymptotic analysis for $\epsilon $ sufficiently small.

Consider the stable leaf of the invariant foliation for (\ref{rde}) ($0<\epsilon \ll 1$), passing through a point $\Phi^0 \in H$,
\begin{align}
W (\Phi^0, \omega) = \{\xi+l^s (\xi,\Phi^0,\omega)\Big| \xi\in H^s\}.       \label{E.mfds}
\end{align}

Let the deterministic leaf (i.e. $\epsilon=0$) be represented as
\begin{align}
 \{\xi+l^{(d)}(\xi)\Big| \xi\in H^s\}, \label{E.mfdd}
\end{align}
where $l^{s}(\cdot,\omega):H^s\rightarrow H^u$ and $l^{(d)}(\cdot) :H^s\rightarrow H^u$ are Lipschitz mappings. We expand
\begin{align}\label{exp_h}
\hse &= l^{(d)}(\xi) +\epsilon \hsk{1} + \epsilon^2 \hsk{2} +\cdots+\epsilon^k \hsk{k}+ \cdots
\end{align}

%Since $\hse$ is determined by (), which can also be written as
%\begin{align}
%h^s(\xi,\omega)&=&\int_{\infty}^0 e^{-As+\int_s^0 z(\theta_r\omega)dr}P^u G(\theta_s\omega,u(s,\xi+h^s(\xi,\omega),\omega))ds\\
%\end{align}
%where $u(t,\xi+h^s(\xi,\omega)$ represents the solution of () with initial condition $u(0)=\xi+h^s(\xi,\omega)$.

Write the solution of (\ref{star_1}) in the form
\begin{align}\label{exp_u}
\psi(t)&=\psi^{(d)}(t)+\epsilon \psi^{(1)}(t)+\cdots+ \epsilon^k \psi^{(k)}(t)+\cdots
\end{align}
with the initial condition
\begin{align}\label{exp_u0}
\psi(0)=\xi+l^s(\xi,\omega)-\Phi^0=\xi -\Phi^0 +l^{(d)}(\xi)+\epsilon l^{(1)}(\xi,\omega)+\cdots.
\end{align}
Similarly, write $\Phi(t)=\phi(t,\omega, \Phi^0)$ as
\begin{align}\label{exp_u1}
\Phi(t)&=\Phi^{(d)}(t)+\epsilon \Phi^{(1)}(t)+\cdots+ \epsilon^k \Phi^{(k}(t)+\cdots
\end{align}
with initial condition
\begin{align}
\Phi(0)=\Phi^0.
\end{align}
By  the Taylor expansion, we obtain
\begin{align}
e^{\zz{t}} =e^{\epsilon Z(\theta_t (\omega))}
  =1+\epsilon \ZZ{t} + \cdots+\frac{\epsilon^{k} \lb\ZZ{t}\rb^{k}}{k!} + \cdots\label{exp_z}
\end{align}
and
\begin{align}
&e^{\sz} =e^{\epsilon \sZ}\notag\\
& \quad =1+\epsilon \sZ + \cdots+\frac{\epsilon^{k} \lb\sZ\rb^{k}}{k!}+\cdots. \label{exp_sz}
 % + \frac{\epsilon^{k+1} \lb\sZ\rb^{k+1}}{(k+1)!} e^{\lambda \sZ}\notag
\end{align}
Suppose $F(u)$ is sufficiently smooth with respect to $u$. With (\ref{exp_z}), it follows from (\ref{E.dfnG}) that
\begin{align}
&G(\theta_t\omega,\psi(t)) =e^{-\zz{t}} F\lb e^{ \zz{t}}\psi(t)\rb\notag\\
&=e^{-\epsilon\ZZ{t}} F\lb e^{ \epsilon\ZZ{t}}\psi(t)\rb\notag\\
&=\lb 1-\epsilon \ZZ{t} +  \cdots \rb F \big( \lb 1+\epsilon \ZZ{t}  +\cdots \rb \lb \psi^{(d)}(t)+\epsilon \psi^{(1)}(t)+\cdots \rb \big)\notag\\
&=F(\psi^{(d)}(t))+\epsilon \lb -\ZZ{t} F(\psi^{(d)}(t))+F_u^{\psi^{(d)}}\lb \psi^{(1)}(t)+\ZZ{t} \psi^{(d)}(t)\rb\rb+\cdots\label{exp_G},
\end{align}
where $F_u^{\psi^{(d)}(t)}$ represents the first order  Fr\'echet derivative \cite{Hunter} of the function $F(u)$ with respect to $u$ and evaluated at $\psi^{(d)}(t)$. In Euclidean space,  the Fr\'echet derivative reduces to the classical derivative.

Substituting (\ref{exp_h}), (\ref{exp_G}) and (\ref{exp_u}) into (\ref{rde}), and equating the terms with the same power of $\epsilon$, we get
\begin{align}
\begin{cases}
\frac{d\psi^{(d)}(t)}{dt} =A\psi^{(d)}(t) +F(\psi^{(d)}(t)+ \Phi^{(d)}(t)-F(\Phi^{(d)}(t)\label{E.u0},\\
\psi^{(d)}(0) =\xi+l^{(d)}(\xi)-\Phi^0,
\end{cases}
\end{align}
and
\begin{align}
\begin{cases}
\frac{d\psi^{(1)}(t)}{dt} =\lB A+F_u^{\psi^{(d)}+\Phi^{(d)}} \rB \psi^{(1)}(t) +\tilde \lambda,  \label{E.u1}\\
\psi^{(1)}(0) =l^{(1)}(\xi,\omega),
\end{cases}
\end{align}
where
\begin{align}
\tilde \lambda &= \ZZ{t} \lB  \psi^{(d)}(t)+F(\Phi^{(d)}(t))-F(\psi^{(d)}(t)+\Phi^{(d)}(t))-\fu{t} \psi^{(d)}(t)\rB \notag\\
&+ F_u^{\psi^{(d)} + \Phi^{(d)}}\lb \Phi^{(1)} (t)+Z(\theta_t\omega) (\psi^{(d)}(t)+\Phi^{(d)}(t))\rb\notag\\
&-F_u^{ \Phi^{(d)}}\lb \Phi^{(1)} (t)+Z(\theta_t\omega) \Phi^{(d)}(t)\rb.
\end{align}

 Solve for $\psi^{(d)}(t)$ and $\psi^{(1)}(t)$,
\begin{align}
\psi^{(d)}(t)=e^{At} \psi^{(d)}(0)+\int_0^t e^{A(t-s)}\lb F(\psi^{(d)}(t)+ \Phi^{(d)}(t)-F(\Phi^{(d)}(t)\rb\,ds,
\end{align}
\begin{align}
\psi^{(1)}(t)=e^{ A t+ \fuss{s}{t} }\lb h^{(1)}(\xi,\omega)-\int_0^t e^{ -A s+ \fus{r}{s} } \tilde \lambda \,ds\rb.
\end{align}
\normalsize

Similarly, we have

\begin{align}
\begin{cases}
\frac{d\Phi^{(d)}(t)}{dt} =A\Phi^{(d)} (t) +F(\Phi^{(d)}(t))\label{E.u011},\\
\Phi^{(d)}(0) =\Phi^0,
\end{cases}
\end{align}
and
\begin{align}
\begin{cases}
\frac{d\Phi^{(1)}(t)}{dt} =\lB A+F_u^{\Phi^{(d)}} \rB \Phi^{(1)}(t)+ \tilde B \label{E.u111}\\
\Phi^{(1)}(0) = 0,
\end{cases}
\end{align}
where
\begin{align}
\tilde B=-\ZZ{t} \lB -\Phi^{(d)}(t)+F(\Phi^{(d)}(t))-F_u^{\Phi^{(d)}} \Phi^{(d)}(t)\rB.
\end{align}

 Moreover, solve for $\Phi^{(d)}(t)$ and $\Phi^{(1)}(t)$,
\begin{align}
\Phi^{(d)}(t)=e^{At} \Phi^{(d)}(0)+\int_0^t e^{A(t-s)}F(\Phi^{(1)}(s))\,ds,
\end{align}
\footnotesize
\begin{align}
\Phi^{(1)}(t)=e^{ A t+ \int_0^t F_u^{\Phi^{(d)}(s)} \,ds}\lb l^{(1)}(\xi,\omega)+\int_0^t e^{ -A s+ F_u^{\Phi^{(d)}} } \tilde B\,ds\rb.
\end{align}
\normalsize

 With (\ref{exp_z}), (\ref{exp_sz}) and (\ref{exp_G}), the right hand side of (\ref{star2}) can be written as
\begin{align}\label{exp_rhs}
P^u \Phi^0 + \int_0^{\infty}  e^{-As+\sz} P^u G(\theta_s\omega,u(s)) \,ds=I_0 + \epsilon I_1+R_2,
\end{align}
where $R_2$ represents the remainder term and the other two terms are, \footnotesize
\begin{align}
I_0&= P^u \Phi^0 + \int_\infty^0 e^{-As} P^u  \lB F(\Phi_0(s)+\Psi^{(d)}(s))-F(\Phi_0(s)) \rB\,ds,  \notag\\
I_1&=\int_\infty^0 e^{-As} \left\{\left(\sZ-\ZZ{s}\right)\,\tilde C\, \right\}\,ds,   \notag
\end{align}
\normalsize with
\begin{align}
\tilde C=&F(\psi^{(d)}(s)+\Phi^{(d)}(s))-F(\Phi^{(d)})(s)+F_u^{\psi^{(d)}+\Phi^{(d)}}\notag\\
 &\times\lb \psi^{(1)}(s)+\Phi^{(1)}(s)+
Z(\theta_s\omega)(\psi^{(d)}+\Phi^{(d)}) \rb.
\end{align}

Substituting (\ref{exp_h}) and (\ref{exp_rhs}) into (\ref{star2}), and matching the powers in $\epsilon$, we get
\begin{align}
l^{(d)}(\xi)=I_0=P^u \Phi^0 + \int_\infty^0 e^{-As} P^u  \lB F(\Phi_0(s)+\Psi^{(d)}(s))-F(\Phi_0(s)) \rB\,ds,
\end{align}
and \footnotesize
\begin{align}
 \hsk{1}=\int_\infty^0 e^{-As} \left\{\left(\sZ-\ZZ{s}\right)\,\tilde C\, \right\}\,ds.   \notag
\end{align}
\normalsize

As a summary, we obtain the following result about approximating invariant foliation for the random evolutionary equation \eqref{rde},
including some random ordinary or partial differential equations.

\begin{theorem}[Approximate  invariant foliation for    random evolutionary equations]
Let
  $W(\Phi^0, \omega)=\{\xi+l^s (\xi,\omega, \Phi^0)\Big| \xi\in H^s\}$ represent
 a stable leaf, passing through a point $\Phi^0$,  of the invariant foliation for the random evolutionary equation $\frac{du}{dt}=Au+z(\theta_t \omega)u+G(\theta_t\omega,u)$. Assume that
\begin{itemize}
\item[(i)] $F(u  )$ is twice continuously Fr\'echet differentiable with respect to $u $;
\item[(ii)] For some $\eta$  ($\alpha>\eta>\beta$), the  following gap condition is satisfied
\begin{align}
K\; L_F \; \lb\frac{1}{\eta-\beta} +\frac{1}{\alpha-\eta}\rb<1.  \label{cdnt_eta}
\end{align}
\end{itemize}
Then for $\epsilon$   sufficiently small, the leaf of the random invariant foliation can be approximated as
\begin{align}
W(\Phi^0,\omega)=\{\xi+ l^{(d)}(\xi) +\epsilon \hsk{1} + R_2  \; \Big| \;  \xi\in H^s\},  \notag
\end{align}
where $\|R_2\|\le C(\omega)\epsilon^2$ with $C(\omega)<\infty, \quad a.s. $,
\begin{align}
l^{(d)}(\xi) = P^u \Phi^0 + \int_\infty^0 e^{-As} P^u  \lB F(\Phi_0(s)+\Psi^{(d)}(s))-F(\Phi_0(s)) \rB\,ds,  \label{E.hd}
\end{align}
and \footnotesize
\begin{align}
&\hsk{1}
 = \int_\infty^0 e^{-As} \left\{\left(\sZ-\ZZ{s}\right)P^u F(u_0) \right.\notag\\
& \quad \times \left. \lb F(\psi^{(d)}(s)+\Phi^{(d)}(s))-F(\Phi^{(d)})(s)+F_u^{\psi^{(d)}+\Phi^{(d)}}\lb \psi^{(1)}(s)+\Phi^{(1)}(s)+
Z(\theta_s\omega)(\psi^{(d)}+\Phi^{(d)}) \rb\rb\right\}\,ds. \label{E.h1}
\end{align}
\normalsize

\end{theorem}

\section{Examples}

Let us look at two examples.

\textbf{Example 1}  \\
Consider a SDE system
\begin{align}\label{example1_1}
\begin{cases}
\dot X=-X+\epsilon X\circ \dot W,\\
\dot Y=Y+X^2+\epsilon Y \circ \dot W,
\end{cases}
\end{align}
where $\epsilon >0$ and $W$ is a scalar Brownian motion.  In this example,
$A=\begin{pmatrix} \phantom{-}-1&\phantom{-}0\\
\phantom{-}0&\phantom{-}1\end{pmatrix}$, $H=\R^2$ (a finite dimensional Hilbert space), $H^s=\left\{\begin{pmatrix}\bar x \\
0\end{pmatrix} \Big| \bar x\in \R\right\}$, and $H^u=\left\{\begin{pmatrix}0 \\
\bar y\end{pmatrix} \Big| \bar y\in \R\right\}$. The transformed differential equations with random coefficients are
\begin{align}\label{ttmmpp}
\begin{cases}
\dot x=-x+\epsilon  Z(\theta_t \omega)x,\\
\dot y=y +\epsilon Z(\theta_t \omega)y +e^{ \epsilon Z(\theta_t \omega)}x^2,
\end{cases}
\end{align}
where $Z(\omega)$ is the stationary solution of $dZ+Zdt=dW$, i.e.~$Z(\omega)~=~\int_{-\infty}^0 e^{\tau} dW(\tau)$ and $Z(\theta_t
\omega)=e^{-t}Z(\omega)+e^{-t}\int_0^t e^\tau dW(\tau)$.

\begin{figure}[htp]
  %\centering
  \epsfig{file=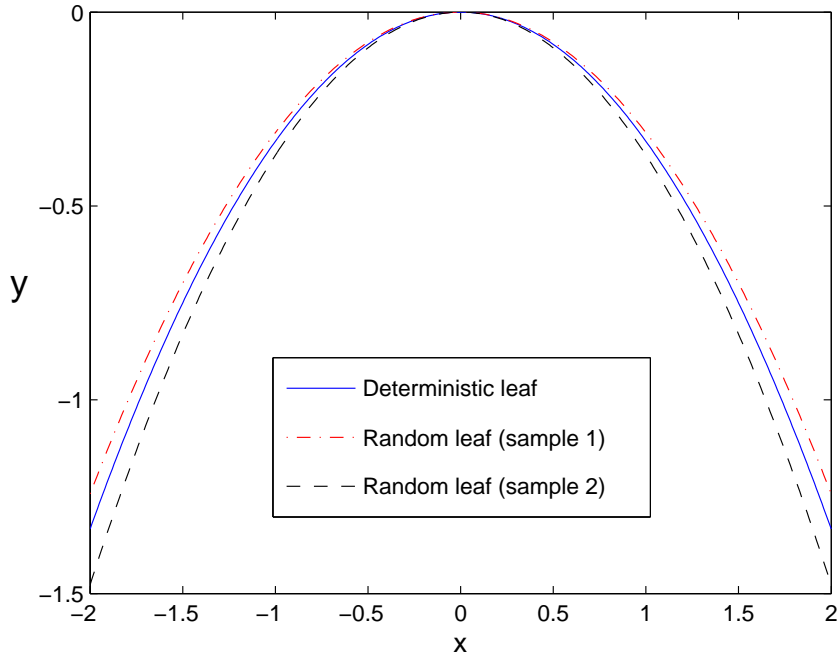,width=\linewidth}
  \caption{The random stable leave passing through the point $(0, 0)$  for Example 1.  Two samples of the random stable leaf are shown here, together with the stable  leaf  for the corresponding deterministic system ($\epsilon=0$)} \label{fig.2}
\end{figure}

The stable leaf of the invariant foliation for (\ref{ttmmpp}), passing through a point $(x_0,y_0)$, can be approximated as \footnotesize
\[
 W = \left\{\begin{pmatrix} x\\  y\end{pmatrix}\Big|   y-y_0 =-\frac{ x^2-x_0^2}{3}-\epsilon  \frac{ x^2-x_0^2}{3} Z(\omega) -\epsilon  \frac{ x^2-x_0^2}{3} \lb\int_0^\infty \ e^{-3
\tau}\,dW_\tau \rb +O(\epsilon^2),\;\; x\in \R\right\}
\]
\footnotesize
%and the corresponding leaf for (\ref{ttmmpp}) can be approximated as \footnotesize
%\begin{eqnarray*}
%& &  W(x_0, y_0, \omega)  =     \nonumber  \\
%& &  \left\{\begin{pmatrix}  x\\  y\end{pmatrix}\Big| y-y_0=-\frac{  x^2-x_0^2}{3}+\epsilon Z(\omega)\lb y_0+\frac{2}{3} x_0^2\rb
%   -    \epsilon
%\frac{ x^2-x_0^2}{3} \lb\int_0^\infty \ e^{-3 \tau}\,dW_\tau \rb
%+O(\epsilon^2)  \;\; x \in \R\right\}.
%\end{eqnarray*}
\normalsize Figure \ref{fig.2} compares a random stable leaf (two samples   are shown here) for (\ref{ttmmpp}) with the deterministic stable leaf,  passing through the point  $x_0=0$
and $y_0=0$.

\textbf{Example 2}  \\
Consider the following SPDE
\begin{align}\notag
\begin{cases}
U_t = \lb U_{xx} +10 U\rb - U^3 +\epsilon U\circ \dot W, x\in [0,1],  \\
U(0,t) =U(1,t)=0,
\end{cases}
\end{align}
where $\epsilon >0$ and $W$ is a scalar Brownian motion. In this example, $A~=~\Delta+10 $, $H=L^2(0,1)$, $D(A)=H^2_0(0,1)$, $F(u )=u^3 $. Note that
the eigenvalues of $A$ are $\lambda_n =10- (n\pi)^2$, and the corresponding normalized eigenfunctions are $e_n=\sqrt{2} \sin(n\pi x)$,
$n=1,\;2\;\cdots$. Here $H^s=Span\left\{e_1\right\}$ and $H^u=Span~\left\{e_2, e_3, \cdots, e_n, \cdots \right\}$.

The transformed random partial differential equation is
\begin{align}\notag
\begin{cases}
u_t = \lb u_{xx} +10 u\rb +Z(\theta_t \omega)u- e^{\epsilon 2 Z(\theta_t\omega)}u^3 , x\in [0,1],  \\
u(0,t) =u(1,t)=0.
\end{cases}
\end{align}

 Unlike Example 1, here we can not express $\hsk{1}$ analytically, but only estimate it via
(\ref{E.h1}).

\end{document}